\documentclass[a4paper]{article}     % Basic GT/GTM/AGT style
\usepackage{amssymb,amsmath,hyperref}
\usepackage{graphics,tikz, enumitem}
\usepackage[english]{babel}

\usepackage{fullpage}

\usepackage[all]{xy}
\title{Epimorphisms of pseudo-quadratic polar spaces}

\author{Petra Schwer \thanks{The first author is supported by the DFG-project SCHW 1550/2-1} \and Koen Struyve \thanks{The second author is supported by the Fund for Scientific Research --
Flanders (FWO - Vlaanderen)}}

\parskip\medskipamount
\newtheorem{thm}{Theorem}[section]
\newtheorem{Theorem}{Main Result}
\newtheorem{rem}[thm]{Remark}
\newtheorem{prop}[thm]{Proposition}
\newtheorem{lemma}[thm]{Lemma}

\def\<{\langle}
\def\>{\rangle}
\newcommand{\proof}{\emph{Proof.~}}

\def\qed{{\hfill\hphantom{.}\nobreak\hfill$\Box$}}

\newcommand{\id}{\mathrm{id}}

\newcommand{\op}{\mathsf{op}}

\begin{document}

\maketitle

\begin{abstract}    % type your abstract below
We classify the epimorphisms of the buildings $\mathsf{BC}_l(K,K_0,\sigma,L, q_0)$, $l \geq 2$, of pseudo-quadratic form type. This completes the %handles the final case of a 
classification of epimorphisms of irreducible spherical Moufang buildings of rank at least two.
\end{abstract}

{\footnotesize AMS Classification: 20E42, 51E12, 51E24}
%\tableofcontents

\section{Introduction}
The aim of this paper is to complete the classification of epimorphisms of irreducible spherical Moufang buildings of rank at least two. For projective planes and spaces defined over a skew field or octonion division algebra $K$ such a classification is given by the work of Andr\'e~\cite{And:69}, Faulkner and Ferrar~\cite{Fau-Fer:83} and Skornjakov~\cite{Sko:57}. It is shown there that such epimorphisms essentially correspond with the total subrings of $K$, i.e. subrings $R \subset K$ such that $K = R \cup (R \setminus \{0\})^{-1}$.  In~\cite{Str:*} the second author derives some of the structure theory of epimorphims of irreducible spherical Moufang buildings of rank at least two and uses this to show that when such a building is defined over a field (for a suitable definition), then these epimorphisms are closely related with affine buildings (and their non-discrete generalizations). 

In the view of these results, the only non-treated case is that of the buildings $\mathsf{BC}_l(K,K_0,\sigma,L, q_0)$ ($l \geq 2$) of pseudo-quadratic form type. The main difference with the cases handled in~\cite{Str:*} is that a total subring of a field always corresponds to a valuation of this field, while this is not true for skew fields in general. As a consequence one can no longer apply the rich theory of affine buildings, meaning that we have to construct the epimorphisms in a different, more ad hoc manner.%, which consists of a large part of this paper.

The precise statement of our classification can be found in Section~\ref{section:main}. 

Finally we note that in this paper we only consider type-preserving epimorphisms between (thick) buildings.

%\section{Preliminaries}
\section{Polar spaces of pseudo-quadratic form type}\label{section:psf}

%\subsection{Polar spaces of pseudo-quadratic form type}\label{section:psf}% TO DO: add another reference to deal with the rank > 2 case?

In this section we define the polar spaces of interest in this paper. Our approach is based on~\cite[(16.5)]{Tit-Wei:02}. Let $K$ be a skew field and $\sigma$ an involution of $K$, meaning $\sigma$ is an anti-automorphism (so $(ab)^\sigma =b^\sigma a^\sigma$) with $\sigma^2 = \id$.  
Let \begin{align*}
&K_\sigma =\{a + a^\sigma \vert a \in K\}, \\
&{K}^\sigma =\{a \vert a \in K, a^\sigma = a\}.
\end{align*}
Choose a $K_\sigma \subset K_0 \subset K^\sigma$ containing the element 1, such that for all $t \in K$ we have $t^\sigma K_0 t = K_0$. Such a set is called an \emph{involutory set}. If the characteristic of $K$ is different from 2, then $K_\sigma = K_0 = K^\sigma$. Let $L$ be a right vector space over $K$.
A map $f:L \times L \rightarrow K$ is a \emph{skew-hermitian sesquilinear form} on $L$ with respect to $\sigma$, if $f(a,b)^\sigma = - f(b,a)$ and $f(at,bu) = t^\sigma f(a,b) u$ for all $a,b \in L$ and $t,u \in K$.  
A map $q: L \rightarrow K$ is a \emph{skew-hermitian pseudo-quadratic form} on $L$ with respect to $\sigma$ if $f$ on $L$ is a skew-hermitian sesquilinear form with respect to $\sigma$, such that the following two conditions are satisfied for all $a,b \in L$ and $t\in K$:

\begin{itemize}
\item $q(a + b) \equiv q(a) +q(b) +f(a,b) \mod K_0$,
\item $q(at) \equiv t^\sigma q(a) t \mod K_0$.
\end{itemize}

If one moreover has that $q(a) \in K_0$ only if $a=0$, then we say that $q$ is \emph{anisotropic}. If all of this is satisfied we say that the quintuple $(K, K_0, \sigma,L,q)$ is an \emph{anisotropic skew-hermitian pseudo-quadratic space}.

\begin{rem} \rm
We will often omit the adjective ``skew-hermitian'' as we will not take other pseudo-quadratic spaces in consideration.
\end{rem}

We are now able to define the rank $l$ polar space $\mathsf{BC}_l(K,K_0,\sigma,L, q_0)$ where $l \geq 2$ is an integer, and $(K,K_0, \sigma,L_0, q_0)$ an anisotropic pseudo-quadratic space. Let $X$ denote the right vector space $L_0 \oplus K^{2l}$. The function $$q: (v \vert a_1, \dots, a_{2l}) \mapsto q_0(v) + a_1^\sigma a_2 + \dots + a_{2l-1}^\sigma a_{2l}$$ with $(v  \vert a_1, \dots, a_{2l}) \in X$, is a pseudo-quadratic form on $X$.
The associated skew-hermitian $f$ is defined by $f((v  \vert a_1, \dots, a_{2l}),(w  \vert b_1, \dots, b_{2l})) = f_0(v,w)$. 

A subspace $S$ of the vector space $X$ is \emph{singular} if $x \in S$ implies $q(x) \in  K_0$. The polar space is now formed by the set of singular subspaces. The points will be the one-dimensional subspaces, the lines the two-dimensional subspace, etc. The building of type $\mathsf{C}_l$ associated to this polar space is the flag complex of this incidence geometry.

%\begin{rem} \rm
%This is a slightly different approach that the one used in~\cite[\S 11]{Tit-Wei:02} (which uses skew-hermitian forms). The advantage of the current approach is that the Main results will also hold for those polar spaces arising from quadratic forms (when $\sigma$ is the identity automorphism).
%\end{rem}

% \subsection{Total subrings}
% 
% A \emph{total subring} of a skew field $K$ is a subring $R$ of $K$ such that $K = R \cup (R \setminus \{0\})^{-1}$. 
% 
% Let $m$ be the set of non-unit elements in such a total subring $R$, then $m$ is the unique maximal (two-sided) ideal of $R$ (see for example~\cite[\S 2]{Fau-Fer:83}). A direct corollary is that the quotient $R  / m$ is a skew field. We call this the \emph{residue skew field} and denote it by $K_R$. 

\section{Statement of the main result}\label{section:main}

We will see in the main theorem below that the total subrings essentially determine the type-preserving epimorphisms of buildings of pseudo-quadratic form type.  Here a \emph{total subring} of a skew field $K$ is a subring $R$ of $K$ such that $K = R \cup (R \setminus \{0\})^{-1}$. 

Let $m$ be the set of non-unit elements in a total subring $R$, then $m$ is the unique maximal (two-sided) ideal of $R$ (see for example~\cite[\S 2]{Fau-Fer:83}). A direct corollary is that the quotient $R  / m$ is a skew field. We call this the \emph{residue skew field} and denote it by $K_R$. 

We will now state the main result.

\begin{Theorem}
Let $(K, K_0, \sigma, L_0, q_0)$ be an anisotropic skew-hermitian pseudo-quadratic space. Every type-preserving epimorphism of the building $\mathsf{BC}_l(K, K_0, \sigma, L_0, q_0)$, $l \geq 2$, is completely determined (up to isomorphisms) by a total subring $R$ of the skew field $K$ and a left coset of $R^*$ in the multiplicative group $K^*$ satisfying the following three conditions.
\begin{enumerate}[label={(C*)}, leftmargin=*]
%\begin{itemize} 
\item[(C1)] The anti-automorphism $a \mapsto a^{\sigma s}$ of $K$ stabilizes $R$,
\item[(C2)]  $(u,t), (w,r) \in T : t,r \in sR \Rightarrow f_0(u,w) \in sR$, 
%\item[(C3)]  $(u,t), (w,r) \in T : t,r \in sm \Rightarrow f_0(u,w) \in sm$,
\item[(C3)]  $(u,t), (w,r) \in T : t \in sR, r\in sm \Rightarrow f_0(u,w) \in sm$,
%\end{itemize}
\end{enumerate}
where  $s$ is an element of the left coset of $R^*$, $f_0$ the skew-hermitian form associated to $q_0$ and $m$ the unique maximal two-sided ideal of $R$. Conversely if $R$ is a total subring of the skew field $K$ and  $sR^*$ a left coset of $R^*$ satisfying these four conditions, then there exists a type-preserving epimorphism of the polar space $\mathsf{BC}_l(K, K_0, \sigma, L_0, q_0)$ determined exactly by this total subring and left coset.
\end{Theorem}

The proof is split into two parts as described in Section~\ref{section:proof}.

\section{Auxiliary results}\label{section:aux}
This section gathers helpful results on spherical buildings, epimorphisms and polar spaces. 

We start by giving root group sequences which describe the rank two Moufang spherical buildings, which are also known as \emph{Moufang polygons}, which appear in buildings of pseudo-quadratic form type. Root group labelings then describe these buildings $\mathsf{BC}_l(K,K_0,\sigma,L, q_0)$ ($l \geq 2$).
In Section~\ref{section:rgl2} we show how the direct construction in Section~\ref{section:psf} relates to the root group labeling.
In Section~\ref{section:sum} we summarize the results from~\cite{Str:*} used in this paper. In particular we describe the interplay between epimorphisms and root group labelings.

\subsection{The root group sequence of $\mathsf{A}_2(K)$}\label{section:a2}
Let $K$ be a skew field. Let $U_i$ ($i\in \{1,2,3\}$) be groups parametrized by isomorphisms $x_i$ from the additive group of $K$ to $U_i$. The only non-trivial commutator relation is given by $$ [x_1(s), x_3(t)] =x_2(st) $$ for $s,t \in K$. This defines a root group sequence $\Theta_{\mathsf{A}_2(K)}$ (see~\cite[(16.1)]{Tit-Wei:02}).

We also list the following identity (from~\cite[(32.5)]{Tit-Wei:02}) which one will need in order to apply Lemma~\ref{lemma:6.7}.
\begin{equation}\label{eq:mu_a}
x_2(u)^{\mu(x_1(t))} = x_3(t^{-1} u) 
\end{equation}

In what follows we will work with $\mathsf{A}_2(K^\op)$. The \emph{opposite} skew field $K^\op$ is defined as the field with the same underlying set as $K$ but with multiplication given by $a * b = ba$ (with $a,b \in K$).

\subsection{The root group sequence of $\mathsf{BC}_2(K,K_0,\sigma,L_0,q)$}\label{section:bc2}
We use the notations from Section~\ref{section:psf}. Let $(K,K_0, \sigma,L_0, q_0)$ be an anisotropic pseudo-quadratic space. Let $T$ be the elements $(w, t)$ in $L_0 \times K$ such that $q_0(w) -t \in K_0$. One derives that if $(w,t) \in T$ then $ f(w,w) = t - t^\sigma $. The set $T$ can be made into a group with multiplication $(w,t) \cdot (v,r) = (w+v, t+r + f(v,w))$ and inverse $(w,t)^{-1} = (-w,-t^\sigma)$. For proofs see~\cite[(11.24), (11.19)]{Tit-Wei:02}.

Let $U_i$ ($i\in \{1,2,3,4\}$) be groups parametrized by the group $T$ in case $i$ is odd and by the additive group of $K$ in case $i$ is even, both via maps $x_i$.  The non-trivial commutator relations are given by:
\begin{align*}
[x_1(w,t), x_3(v,r)^{-1}] &= x_2(f_0(w,v)), \\
[x_2(k), x_4(a)^{-1}] &= x_3(0, k^\sigma a + a^\sigma k), \\
[x_1(w,t), x_4(k)^{-1}] &= x_2(tk) x_3(wk,k^\sigma t k),  
\end{align*}
for $(w,t), (v,r) \in T$ and $k,a \in K$. These relations define a root group sequence $\Theta_{\mathsf{BC}_2(K,K_0,\sigma,L_0,q)}$.

We end by giving the following equations from~\cite[(32.9)]{Tit-Wei:02}.
\begin{align}
x_1(w,t)^{\mu(x_4(k))} =& x_3(-wk, k^\sigma t k) \label{eq:mu3_bc} \\
x_2(k)^{\mu(x_4(a))} =& x_2(- a^{-\sigma} k^\sigma a ) \label{eq:mu2_bc} \\
x_4(a)^{\mu(x_1(w,t))} =& x_2(t a)  \label{eq:mu_bc} 
\end{align}

\subsection{The root group labeling of $\mathsf{BC}_l(K,K_0,\sigma,L, q_0)$}\label{section:rgl}
In this section we describe the root group labeling $(u,\Theta,\theta)$ of the building $\mathsf{BC}_l(K,K_0,\sigma,L, q_0)$, following~\cite[(12.12), (12.16)]{Wei:03}. We will not give every detail of it, only the parts  relevant for our proof. Let $\Pi$ be the following Coxeter diagram with numbered vertices. 

\begin{center}
\begin{tikzpicture}
\fill (0,0) circle (2pt) node[anchor=south] {$1$};
\fill (1,0) circle (2pt) node[anchor=south] {$2$};
\fill (3,0) circle (2pt) node[anchor=south] {$l-2$};
\fill (4,0) circle (2pt) node[anchor=south] {$l-1$};
\fill (5,0) circle (2pt) node[anchor=south] {$l$};

\draw (0,0) -- (1.3 ,0) (2.7,0) -- (4,0) (4,0.05) -- (5,0.05) (4,-.05) -- (5,-.05); 
\draw[style=dotted] (1.5,0) -- (2.5 ,0); 
\end{tikzpicture}
\end{center}

For $i \in \{1,2, \dots, l-1 \}$, let $u(i)$ be isomorphic with the additive group of the skew field $K$. We set $u(l)$ to be isomorphic with the group $T$. We parametrize the groups $u(i)$ ($i \in \{1,2, \dots, l \}$) by isomorphisms $y_i$ from $K$ or $T$ (where applicable) to $u(i)$.

Let $\Theta_{i,i+1} = \Theta_{\mathsf{A}_2(K^\op)}$  for $i\in \{1,2,\dots, l-2\}$ and $\Theta_{l,l- 1} =\Theta_{\mathsf{BC}_2(K,K_0,\sigma,L_0,q)}$. This defines the root group labeling $(u,\Theta,\theta)$ of the building $\mathsf{BC}_l(K,K_0,\sigma,L, q_0)$.

\subsection{Realization of the root group labeling}\label{section:rgl2}
We will now show how the direct construction of $\mathsf{BC}_l(K,K_0,\sigma,L, q_0)$ given in Section~\ref{section:psf} realizes the root group labeling given in Section~\ref{section:rgl} in the sense of~\cite[(12.10-11)]{Wei:03}. We do this by showing how the groups $u(i)$ from the root group labeling act on the vector space $X$.
\begin{align*}
(v \vert a_1, \dots, a_{2l}&)^{y_1(k)} = (v \vert a_1 , \dots , a_{2l-3} +k a_{2l-1} ,a_{2l-2} , a_{2l-1} , a_{2l} -k^\sigma a_{2l-2}), \\
&\cdots \\
(v \vert a_1, \dots, a_{2l}&)^{y_{l-1}(k)} = (v \vert a_1 +k a_3,a_2 , a_3, a_4 -k^\sigma a_2 , \dots , a_{2l}), \\
(v \vert a_1, \dots, a_{2l}&)^{y_l(w,t)} = (v+w a_1 \vert a_1,a_2 -  t a_1  - f_0(w,v),a_3 , \dots , a_{2l}),
\end{align*}
where $k \in K$ and $(w,t) \in T$. The omitted coordinates are left invariant. These maps fix the chamber consisting of the subspaces 
\begin{align*}
\langle (0 \vert &0 , \dots, 0, 1) \rangle, \\
\langle (0 \vert &0 , \dots, 0, 1), (0 \vert 0 , \dots,0,1 ,0, 0) \rangle, \\
\langle (0 \vert &0 , \dots, 0, 1), (0 \vert 0 , \dots, 0,1 ,0, 0), (0 \vert 0 , \dots,0, 1,0,0 ,0, 0) \rangle, \\ 
& \cdots  
\end{align*}

\begin{rem} \rm \label{remark:choice}
Note that the maps of the form 
\begin{align*}
\zeta_i(m): (v \vert \dots, a_{2i-2}, & a_{2i-1},a_{2i}, a_{2i+1}, \dots) \mapsto \\
&(v \vert \dots, a_{2i-2}, m a_{2i-1}, m^{-\sigma} a_{2i}, a_{2i+1}, \dots) 
\end{align*}
with $m \in K^*$ and $i \in \{1,\dots, l\}$ induce automorphisms of the polar space. This map $\zeta_i$ normalizes each of the groups $u(j)$ ($j \in \{1,\dots, l\}$) and acts on the groups $u(j)$ ($j \in \{1,\dots, l-1\}$) as follows.
$$y_j(k)^{\zeta_i(m)} =  \left\{\begin{array}{cl}
y_j(mk) & \mbox {if } j = l-i \\
y_j(km^{-1}) & \mbox {if } j = l-i+1 \\
y_j(k) & \mbox{otherwise} \\
\end{array}\right.$$ 
By combining these automorphisms one can assume without loss of generality that $y_j(1) \in v(j) \setminus w(j)$ for all $j \in \{1,\dots, l-1\}$, where $v(j)$ and $w(j)$ are subgroups of $u(j)$ which will be  introduced in Section~\ref{section:sum}.
\end{rem}

% and correspond with the choice of apartment in $\mathsf{BC}_l(K,K_0,\sigma,L, q_0)$ spanned by the points 
%\begin{align*}
%\langle (0 \vert &1,0 , \dots, 0) \rangle, \\
%\langle (0 \vert &0,1,0 , \dots, 0) \rangle, \\
%\langle (0 \vert &0,0,1,0 , \dots, 0) \rangle, \\
%& \dots 
%\end{align*}

\subsection{A summary of results on epimorphisms of spherical Moufang buildings}\label{section:sum}
In this section we summarize the results in~\cite{Str:*} that we use in the current paper. 

Let $\Delta$, $\Delta'$ be two irreducible spherical Moufang buildings of rank at least two and $\phi$ a (type-preserving) epimorphism from $\Delta$ to $\Delta'$. 

We start by the rank two case. Let $c$ be a chamber in some apartment $\Sigma$ of $\Delta$. With this choice of chamber and apartment there corresponds a root group sequence $(U_+,U_1, \dots, U_n)$. Section 6.1 of~\cite{Str:*} states that the epimorphism $\phi$ induces subgroups $W_i \lhd  V_i \leq U_i$ for every $i$. The subgroup $V_i$ consists of those automorphisms $g \in U_i$ such that there exists an automorphism $g'$ of $\Delta'$ making the following diagram commute.

$$\xymatrix{ \Delta \ar[r]^{g}  \ar[d]_\phi &  \Delta  \ar[d]^\phi  \\
\Delta' \ar[r]_{g'}   &  \Delta'}$$

The subgroup $W_i$ is then the subgroup of those elements in $V_i$ such that the corresponding $g'$ in the previous diagram is the identity automorphism.

The following three lemmas describe how these different subgroups are related.

\begin{lemma}\label{lemma:6.7}
Let $v_i \in V_i\setminus W_i$, then 
\begin{align*}
V_j^{\mu(v_i)} &= V_{2i+n -j}, \\
W_j^{\mu(v_i)} &= W_{2i+n -j}
\end{align*}
for each $i,j \in \{1,\dots,n \}$ such that $2i+n -j \in \{1,\dots,n \}$. 
\end{lemma}
\proof
See~\cite[Cor. 6.7]{Str:*}. 
\qed

\begin{lemma}\label{lemma:6.8}
Choose $u_1 \in U_1$ and $u_n \in V_{n} \setminus W_{n}$.  Let $[u_1,u_n^{-1}] = u_2 \dots u_{n-1}$ (with $u_i \in U_i$), then 
\begin{align*}
u_1 \in V_1 &\Leftrightarrow  u_2 \in V_2, \\
u_1 \in W_1 &\Leftrightarrow  u_2 \in W_2.
\end{align*}
\end{lemma} 
\proof
See~\cite[Lem. 6.8]{Str:*}. 
\qed

\begin{lemma}\label{lemma:6.5}
Choose $u_1 \in V_1$ and $u_3 \in W_3$.  If $[u_1,u_3] = u_2$ then $u_2 \in W_2$.
\end{lemma} 
\proof
This is a special case of~\cite[Cor. 6.5]{Str:*}. 
\qed

The arbitrary rank case can now be approached as follows. Choose a chamber $c$ of the building $\Delta$ and let $(u,\Theta,\theta)$ be a root group labeling associated with this choice of chamber (see~\cite[(12.10-11)]{Wei:03}). The epimorphism $\phi$ again induces subgroups $w(i) \lhd  v(i) \leq u(i)$ for every $i$ as before. These subgroups determine the structure of $\phi$, as shown by the following lemma.

\begin{lemma}\label{lemma:6.12}
If the subgroups $w(i) \lhd v(i)\leq u(i)$ are known for a root group labeling $(u,\Theta,\theta)$ of a spherical Moufang building $\Delta$, then the corresponding epimorphism of $\Delta$ is unique up to isomorphisms.
\end{lemma}
\proof
See~\cite[Cor. 6.12]{Str:*}. 
\qed

For a root group sequence $\Theta_{ij}$ of the root group labeling the $u(i)$ and $u(j)$ form the extremal root groups $U_1$ and $U_n$ of this root group sequence. The subgroups $w(i) \lhd  v(i) \leq u(i)$ and $w(j) \lhd  v(j) \leq u(j)$ correspond respectively to the subgroups $W_i \lhd  V_i \leq U_i$ and $W_j \lhd  V_j \leq U_j$.

\begin{lemma}\label{lemma:proj}
If a certain label $i$ corresponds with a rank one residue which is a projective line over a skew field $K$ with $u(i)$ indexed by $K$ via a map $y_i$, then there exists a total subring $R$ of $K$ with maximal ideal $m$ and a constant $a \in K$  such that 
\begin{align*}
v(i) &= \{y_i(k) | k \in R a \}, \\
w(i) &= \{y_i(k) | k \in m a \}.
\end{align*}
\end{lemma}
\proof
See~\cite[Lem. 7.2-3]{Str:*}.
\qed

\subsection{Some properties of polar spaces}
In this section we state some properties of polar spaces of (pseudo-)quadratic form type needed later on. 

\begin{rem} \rm
We will always suppose in this section that our polar spaces are non-singular and not of hyperbolic type. In particular this is the case for polar spaces corresponding to a building of type $\mathsf{C}_l$.
\end{rem}

Each set of mutually collinear points as well as each subspace of a rank $l$ polar space is contained in a (maximal) subspace of (geometric) dimension $l-1$. These maximal subspaces are called the \emph{generators}. %A subspace of dimension at most $t$ ($t < l$) is contained in at least three subspaces of dimension $t+1$. (One way to derive these properties is by interpreting the flag complex as a building of type $\mathsf{C}_l$ and the fact that subspaces of co-dimension one correspond to vertices in the building.)

\begin{lemma}\label{lemma:fnd}
A subspace of dimension $l-2$ is contained in at least three generators.
\end{lemma}
\proof
This follows from the thickness of the building associated with the polar space.
\qed

The following two lemmas show how points and generators interact.

\begin{lemma}\label{lemma:polar}
Given a generator $\pi$ and a point $p$ not in $\pi$, there is a unique $l$-dimensional subspace $\xi$ containing $p$ and intersecting $\pi$ in a subspace of co-dimension 1. This subspace consists exactly of the points of $\pi$ collinear with $p$.
\end{lemma}
\proof
This property is part of the incidence geometric definition of polar spaces, see for example~\cite[p. 556]{hand}. \qed

\begin{lemma} \label{lemma:polar2}
Let $\pi$ be a $t$-dimensional subspace and $p$ a point  not in this subspace. The set of points in $\pi$ collinear with $p$ either forms a $(t-1)$-dimensional subspace, or every point of $\pi$ is collinear with $p$. Moreover each $(t-1)$-dimensional subspace of $\pi$ arises in this way.
\end{lemma}
\proof
The first assertion follows directly from Lemma~\ref{lemma:polar}. In order to prove the second assertion let $\zeta$ be a $(t-1)$-dimensional subspace of $\pi$ and embed $\pi$ in a generator $\xi$. We then can find a subspace $\chi$ of co-dimension 1 in $\xi$ such that the intersection of $\xi$ and $\chi$ is exactly $\zeta$. Lemma~\ref{lemma:fnd} allows us to find a generator $\xi'$ containing $\chi$ and different from $\xi$. If $p$ is a point of $\xi'$ not in $\chi$, then the points of $\xi$ collinear with $p$ have to be exactly the points of $\chi$ by Lemma~\ref{lemma:polar}. Restricting to the subspace $\pi$ of $\xi$ shows that $\zeta$ consists exactly of those points of $\pi$ collinear with $p$.
\qed

%\begin{lemma}\label{lemma:polar1}
%Given an $n$-dimensional subspace $\alpha$ and a point $p$ in it, there exists a set $P$ of $n-1$ points such that $x$ is the only point of $\alpha$ collinear to all points in $P$.
%\end{lemma}
%\proof
%Embed $\alpha$ in a generator $\pi$. One can choose $n-1$ $(l-1)$-dimensional subspaces of the generator such that the intersection of all of these with $\alpha$ is exactly $p$. For each of these subspaces choose a generator containing it different from $\pi$. Choose a point not on $\pi$ in each of these generators. The previous lemma implies that this set of points has the desired property. \qed
%
%\begin{lemma}\label{lemma:polar2}
%Given an $n$-dimensional subspace $\alpha$ and a set $P$ of $n-1$ points of a polar space. Then there exists a point of $\alpha$ collinear with all points in $P$.
%\end{lemma}
%\proof
%Embed $\alpha$ in a generator $\pi$. The points of $\pi$ collinear with all points of $P$ is the intersection of $n-1$ subspaces of dimension $l-1$ of $\pi$ by Lemma~\ref{lemma:polar}. This intersection is subspace of dimension at least $l-n+1$ of $\pi$, which necessarily intersects the $n$-dimensional subspace $\alpha$ in a one-dimensional space. Hence the lemma is proven. \qed
%

\subsection{Collinearity in $\mathsf{BC}_l(K,K_0,\sigma,L, q_0)$}\label{section:coll}

One checks that the points $ \langle (v\vert a_1,a_2, \dots, a_{2n}) \rangle $ and $\langle (w\vert b_1,b_2, \dots, b_{2n}) \rangle $ of the %polar
space $\mathsf{BC}_l(K,K_0,\sigma,L, q_0)$ are collinear if and only if 
$$f_0(v,w) + a_1^\sigma b_2 + b_1^\sigma a_2 + \dots  + a_{2n-1}^\sigma b_{2n} + b_{2n-1}^\sigma a_{2n} =0.$$
The left-hand side of this equation is the skew-hermitian sesquilinear form associated to the pseudo-quadratic form $q$ on $X$ applied to the two vectors.

\subsection{Polar spaces of quadratic form type}

In this section we define the polar spaces $\mathsf{B}_l(K, L_0, q_0)$ of quadratic form type. We do this as these polar spaces will arise as images of epimorphisms in Section~\ref{section:suf}.

A \emph{quadratic space} $(K,L_0,q_0)$ is a triple consisting of a field $K$, a non-trivial vector space $L_0$ over $K$, equipped with an \emph{quadratic form} $q_0$. This is a map $q:L_0 \to K$ such that there exists a (necessarily unique) bilinear form $f$ on $L_0$ satisfying the following two properties:
\begin{itemize}
\item
$q(u+v) = q(u) + q(v) +f(u,v)$,
\item
$q(tu) = t^2 q(u)$,
\end{itemize}
for all $u,v \in L_0$. The quadratic form $q_0$ is \emph{anisotropic} (and $(K,L_0,q_0)$ an \emph{anisotropic quadratic space}) if one has for every $u \in L_0$ that $q(u) = 0$ if and only if $u=0$.

We can now define the rank $l$ polar space $\mathsf{B}_l(K, L_0, q_0)$ where $ l \geq 2$ is an integer and $(K,L_0,q_0)$ an anisotropic quadratic space. Let $X$ denote the vector space $L_0 \oplus K^{2l}$. The map
$$q : (v \vert a_1, \dots, a_{2l}) \longmapsto q_0(v) +a_1a_2 + \dots + a_{2l-1} a_{2l}$$
is a quadratic form on $X$. A subspace $S$ is called \emph{singular} if it is mapped to zero by $q_0$. As in the pseudo-quadratic form case, the polar space is formed by the singular subspaces and the associated building is the flag complex of the polar space.

\section{Proof of the main result}\label{section:proof}

We split the proof of the main result in two parts. In Section~\ref{section:nec} we will prove:

\begin{thm}\label{thm:main1}
A (type-preserving) epimorphism of a polar space $\mathsf{BC}_l(K, K_0, \sigma, L_0, q_0)$ is completely determined (up to isomorphisms) by a total subring $R$ of the skew field $K$ and a left coset of $R^*$ in the multiplicative group $K^*$ satisfying the following three conditions.
\begin{enumerate}[label={(C*)}, leftmargin=*]
%\begin{itemize} 
\item[(C1)] The anti-automorphism $a \mapsto a^{\sigma s}$ of $K$ stabilizes $R$,
\item[(C2)]  $(u,t), (w,r) \in T : t,r \in sR \Rightarrow f_0(u,w) \in sR$, 
%\item[(C3)]  $(u,t), (w,r) \in T : t,r \in sm \Rightarrow f_0(u,w) \in sm$,
\item[(C3)]  $(u,t), (w,r) \in T : t \in sR, r\in sm \Rightarrow f_0(u,w) \in sm$,
%\end{itemize}
\end{enumerate}
where $s$ is an element of the left coset of $R^*$, $f_0$ the skew-hermitian form associated to $q_0$ and $m$ the unique maximal two-sided ideal of $R$.
\end{thm}

Section~\ref{section:suf} is devoted to the proof of the following theorem.

\begin{thm}\label{thm:main2}
Consider the polar space $\mathsf{BC}_l(K, K_0, \sigma, L_0, q_0)$. If $R$ is a total subring of the skew field $K$ and  $sR^*$ a left coset of $R^*$ in the multiplicative group $K^*$ satisfying the following three conditions.
\begin{enumerate}[label={(C*)}, leftmargin=*]
%\begin{itemize} 
\item[(C1)] The anti-automorphism $a \mapsto a^{\sigma s}$ of $K$ stabilizes $R$,
\item[(C2)]  $(u,t), (w,r) \in T : t,r \in sR \Rightarrow f_0(u,w) \in sR$, 
%\item[(C3)]  $(u,t), (w,r) \in T : t,r \in sm \Rightarrow f_0(u,w) \in sm$,
\item[(C3)]  $(u,t), (w,r) \in T : t \in sR, r\in sm \Rightarrow f_0(u,w) \in sm$,
%\end{itemize}
\end{enumerate}
where $f_0$ is the skew-hermitian form associated to $q_0$ and $m$ the unique maximal two-sided ideal of $R$, then there exists a (type-preserving) epimorphism of the polar space $\mathsf{BC}_l(K, K_0, \sigma, L_0, q_0)$ for which Theorem~\ref{thm:main1} gives rise to the same total subring and left coset.
\end{thm}

Once both of the theorems are proved, the main result follows by combining these.

\subsection{Proof of Theorem~\ref{thm:main1}}\label{section:nec}

Let $\Delta$ be the building $\mathsf{BC}_l(K,K_0,\sigma,L, q_0)$, $f_0$ the skew-hermitian form associated to $q_0$, $(u,\Theta,\theta)$ its root group labeling as given in Section~\ref{section:rgl} and $\phi$ a type-preserving epimorphism from $\Delta$ to another building $\Delta'$ of type $\mathsf{C}_l$.

By Section~\ref{section:sum} we know that this epimorphism is essentially described by subgroups $w(i) \lhd v(i) \leq u(i)$ for $i \in \{1, \dots ,l\}$. Remark~\ref{remark:choice} allows us to assume without loss of generality that $y_i(1) \in v(i) \setminus w(i)$ for $i \in \{1, \dots, l-1\}$. We also fix an element $y_l(v,s) \in v(l) \setminus w(l)$. We pick this element such that $v=0$ if possible.

By Lemma~\ref{lemma:proj} we know that there exists a total subring $R$ of $K$ with maximal ideal $m$ and a constant $a \in K$ such that
\begin{align*}
v(1) &= \{y_1(k) | k \in R a \}, \\
w(1) &= \{y_1(k) | k \in m a \}.
\end{align*}

The next lemma extends these expressions for other $u(i)$, and shows that one can assume that $a =1$. 

\begin{lemma}\label{lemma:n1}
For every $i \in \{1, \dots, l-1 \}$ one has
\begin{align*}
v(i) &= \{y_i(k) | k \in R  \}, \\
w(i) &= \{y_i(k) | k \in m  \}.
\end{align*}
\end{lemma}
\proof
We proof this by induction. We first consider the case $i=1$. As $y_1(1) \in v(1) \setminus w(1)$ it follows that $a^{-1}$ (and so also $a$) is a unit of $R$, and that the statement is true for $i=1$.

Now suppose that the statement is true for some $j \in \{1, \dots, l-2\}$. From Section~\ref{section:rgl} we know that $\Theta_{j,j+1} = \Theta_{\mathsf{A}_2(K^\op)}$. Hence we can identify the subgroups $w(j) \lhd  v(j) \leq u(j), w(j+1) \lhd v(j+1) \leq u(j+1)$ with groups $W_1 \lhd V_1 \leq U_1, W_3 \lhd  V_3 \leq U_3$, respectively, as outlined in Section~\ref{section:sum}, and $U_1$ and $U_3$ as in Section~\ref{section:a2}. These identifications imply that
\begin{align*}
V_1 &= \{x_1(k) | k \in R  \}, \\
W_1 &= \{x_1(k) | k \in m  \}, \\
x_3&(1) \in V_3 \subset W_3. 
\end{align*}
Applying Lemma~\ref{lemma:6.8} and the commutator relation $ [x_1(b), x_3(1)^{-1}] =x_2(-b)$ for $b$ in $K$ we see that
\begin{align*}
V_2 = \{x_2(k) | k \in R  \}, \\
W_2 = \{x_2(k) | k \in m  \}.
\end{align*}
From Equation~(\ref{eq:mu_a}) in Section~\ref{section:a2} we know that $x_2(u)^{\mu(x_1(1))} = x_3(u)$, so Lemma~\ref{lemma:6.7} yields
\begin{align*}
V_3 = \{x_3(k) | k \in R  \}, \\
W_3 = \{x_3(k) | k \in m  \},
\end{align*}
which is, via the identifications, exactly what we need to prove. \qed

The next lemma determines the subgroups $w(l)$ and $v(l)$.

\begin{lemma}\label{lemma:n2}
The subgroups $v(l)$ and $w(l)$ are described by
\begin{align*}
v(l) &= \{y_l(w,t) | (w,t) \in T,  t \in sR  \}, \\
w(l) &= \{y_l(w,t) | (w,t) \in T, t \in s m  \}. 
\end{align*}
\end{lemma}
\proof
As $\Theta_{l,l- 1} =\Theta_{\mathsf{BC}_2(K,K_0,\sigma,L_0,q_0)}$ (see Section~\ref{section:rgl}), one can identify the subgroups $w(l-1) \lhd  v(l-1) \leq u(l-1)$ and $w(l) \lhd  v(l) \leq u(l)$  with groups $W_4 \lhd V_4 \leq U_4$ and $W_1 \lhd  V_1 \leq U_1$, respectively, as outlined in Section~\ref{section:sum}, where $U_1$ and $U_4$ are as in Section~\ref{section:bc2}.

Lemma~\ref{lemma:n1}  implies that $V_4$ and $W_4$ can be expressed as 
\begin{align*}
V_4 = \{x_4(k) | k \in R  \}, \\
W_4 = \{x_4(k) | k \in m  \}. 
\end{align*}

By Lemma~\ref{lemma:6.7}, Equation~(\ref{eq:mu_bc}) (see Section~\ref{section:bc2}) and $y_l(v,s) \in v(l) \setminus w(l)$ one obtains that 
\begin{align*}
V_2 = \{x_2(k) | k \in sR  \}, \\
W_2 = \{x_2(k) | k \in s m  \}. 
\end{align*}
It is now possible to describe the relevant subgroups of $U_1$ using Lemma~\ref{lemma:6.7} and the commutator relation $[x_1(w,t), x_4(1)^{-1}] = x_2(t) x_3(w, t )$ found in Section~\ref{section:bc2}. One derives that $x_1(w,t) \in V_1$ or $W_1$ if and only if $t \in sR$ or $sm$ respectively, so
\begin{align*}
V_1 = \{x_1(w,t) | (w,t) \in T, t \in sR  \}, \\
W_1 = \{x_1(w,t) | (w,t) \in T, t \in s m  \},
\end{align*}
which is what we need to show. \qed

At this point we have determined all of the subgroups $w(i) \lhd v(i) \leq u(i)$ for $i \in \{1, \dots, l\}$. These subgroups are completely encoded by the total subring $R$ and an element  $s \in K^*$ (or more exactly a left coset of $R^*$ in the multiplicative group $K^*$). These subgroups determine on their turn the epimorphism by Lemma~\ref{lemma:6.12}.  

In the remainder of this section we will derive the properties that these $R$ and $s$ satisfy.  

\begin{lemma}\label{lemma:stab}
The map $a \mapsto s^{-1} a^{\sigma} s^\sigma$ of $K$ stabilizes $R^*$.
\end{lemma}
\proof
We use the same setting of Lemma~\ref{lemma:n2}. From Lemma~\ref{lemma:6.7} and $x_2(s)^{\mu(x_4(a))} = x_2(- a^{-\sigma} s^\sigma a )$ (see Equation~(\ref{eq:mu2_bc})), where we pick $a \in R^*$ and $k$ equal to $s$, it follows that $-a^{-\sigma} s^\sigma a \in sR $. As $a$ is an invertible element of $R$ and $a^{-1} \in R^*$ this is equivalent to $s^{-1} a^{\sigma} s^\sigma \in R$. \qed

\begin{prop}\label{prop:c1}
The anti-automorphism $a \mapsto a^{\sigma s}$ of $K$ stabilizes $R$.
\end{prop}
\proof
As this map is the combination of an automorphism and anti-automorphism it is clear that it is an anti-automorphism. Let $a \in R^*$. Then $a^{\sigma s} = s^{-1} a^\sigma s = (s^{-1} a^\sigma s^\sigma) (s^{-\sigma} s)$. The first factor and the inverse of the second factor are of the form as in Lemma~\ref{lemma:stab}, so both of them and their product lie in $R^*$. So the anti-automorphism $a \mapsto a^{\sigma s}$ maps $R^*$ into $R^*$. As $R^*$ generates  $R$ as a ring it follows that $R$ is mapped into $R$.

The inverse of the map $a \mapsto a^{\sigma s}$ is given by the map $a \mapsto s^{-\sigma}a^{\sigma } s^\sigma$. One shows analogously, using the decomposition $s^{-\sigma}a^{\sigma } s^\sigma = (s^{-1}s^{\sigma} ) s^{-\sigma}a^{\sigma } s^\sigma$, that this inverse maps $R$ into $R$. Hence we conclude that the anti-automorphism $a \mapsto a^{\sigma s}$ of $K$ stabilizes $R$.
\qed

%One checks easily that the map $a \mapsto s^{-1} a^{\sigma} s^\sigma$ is involutory, hence it stabilizes 

\begin{prop}\label{prop:c2}
$$ \forall (u,t), (w,r) \in T : t,r \in sR \Rightarrow f_0(u,w) \in sR. $$
%\forall (u,t), (w,r) \in T : t,r \in sm \Rightarrow f_0(u,w) \in sm.
%\end{align*}
\end{prop}
\proof 
Let $ (u,t), (w,r) \in T$ such that  $t,s \in sR$. Lemma~\ref{lemma:n2} implies that $y_l(u,t) , y_l(w,r) \in v(l)$. As $v(l)$ is a subgroup it follows that the product $y_l(w,r) \cdot  y_l(u,t)$ also lies in $v(l)$, hence $t + r + f_0(u,w) \in sR$ (see Section~\ref{section:bc2} and again Lemma~\ref{lemma:n2}). Because $R$ is a ring this is equivalent with $f_0(u,w) \in sR$.  \qed

\begin{prop}\label{prop:c3}
$$\forall (u,t), (w,r) \in T : t \in sR, r\in sm \Rightarrow f_0(u,w) \in sm. $$
\end{prop}
\proof
We use the same setting of Lemma~\ref{lemma:n2}. Let $(u,t), (w,r) \in T$ with $t\in sR, r \in sm$. Note that $x_1(u,t) \in V_1$ and $x_1(w,r) \in W_1$. We start by determining $V_3$ using Lemma~\ref{lemma:6.7}, Equation~(\ref{eq:mu3_bc}) and $x_4(1) \in V_4 \setminus W_4$. Combining this yields that $x_3(w,r) \in W_3$. Lemma~\ref{lemma:6.5} now implies that $[x_1(u,t), x_3(w,r)^{-1}] = x_2(f_0(u,w))  \in W_2$ which is equivalent to $f_0(u,w) \in sm$. \qed

As Propositions~\ref{prop:c1},~\ref{prop:c2} and~\ref{prop:c3} prove Conditions (C1)-(C3), this concludes the proof of Theorem~\ref{thm:main1}.

\subsection{Proof of Theorem~\ref{thm:main2}}\label{section:suf}

In this section we construct epimorphisms of the polar space $\mathsf{BC}_l(K, K_0, \sigma, L_0, q_0)$. We do this starting from a total subring $R \subset K$ and a left coset $sR^*$ of $R^*$ in the multiplicative group $K^*$ satisfying the conditions outlined in Theorem~\ref{thm:main2}.

As before we will let $m$ denote the set of non-invertible elements of $R$ and $K_R$ the corresponding residue skew field. 

%\begin{itemize
%The next step in the proof is to start with data satisfying the conditions derived in Theorem~\ref{thm:main1} and construct an epimorphism of the polar space $\mathsf{BC}_l(K, K_0, \sigma, L_0, q_0)$ given a total subring $R \subset K$ satisfying the conditions derived in the previous section. A %We assume that the compatibility conditions derived in the previous section are satisfied.

%Note that we can suppose that $l=1$ by Remark~\ref{remark:choice}.

\subsubsection{Structure of $K_R$}

We start by showing that one can choose the representative $s$ in the left coset in a special way.

\begin{lemma}\label{lemma:cases}
The left coset $sR^*$ contains an element $r$ such that we are in exactly one of the following two cases:
\begin{enumerate}[label={(Case **.)}, leftmargin=*]
%\begin{itemize}
\item[Case I.]
$r \in K_0$ and $a \mapsto a^{\sigma r}$ is an involution,
\item[Case II.]
$ K_0 \cap rR^*  = \varnothing$, $r^{-1}r^\sigma + 1 \in m$, $K_R$ is a field and $a \mapsto a^{\sigma r}$ induces the identity on $K_R$. 
%\end{itemize}
\end{enumerate}
\end{lemma}
\proof
Note that any element $r$ in $K_0 \cap sR^* $ is fixed by $\sigma$, implying that $a \mapsto a^{\sigma r}$ is an involution. So such an element directly satisfies Case I. Hence we may suppose that $K_0 \cap sR^*$ is empty. This implies that we have two possibilities for $1+s^{-1}s^\sigma =  s^{-1}(s+s^\sigma) \in s^{-1}K_\sigma \subset s^{-1}K_0$. It can either be an element of $m$, or an element of $K \setminus R$. Suppose the latter holds. Then $s^{-1}s^\sigma$ also belongs to $k \setminus R$, and the inverse $s^{-\sigma} s$ belongs to $m$. But $(s^{-1}s^\sigma)^{\sigma s} = s^{- \sigma} s$ which contradicts Condition (C1). So $1+s^{-1}s^\sigma \in m$ and consequently $s^{-1}s^\sigma \in R$.  Now for a given $a \in R$ we have that $s^{-1}(s^\sigma a + a^\sigma s ) = s^{-1} s^\sigma a + a^{\sigma s} \in s^{-1}K_0 $. This element lies in $m$ as it can not be an element of $k\setminus R$ (by Condition (C1) and what we already had derived for $s^{-1}s^\sigma$). Combined with $1+s^{-1}s^\sigma \in m$ this implies that $a \equiv a^{\sigma s} \mod m$, or that $a \mapsto a^{\sigma r}$ is the identity automorphism of $K_R$. Because it is also an anti-automorphism this yields that $K_R$ is a field. Finally we remark that  $r \in K_0$ and $K_0 \cap rR^*  = \varnothing$ are mutually exclusive, so exactly one of the two cases holds.
\qed

From now on suppose that $s$ is as in one of the two cases described by this lemma. We denote the anti-automorphism induced on $K_R$  by the anti-automorphism $a \mapsto a^{\sigma s}$ by $\sigma_R$.

\begin{lemma}
Under the assumption of Case I,  $\overline{K_0} := s^{-1}K_0 \cap R \mod m$ is an involutory set of $K_R$.
\end{lemma}
\proof 
One calculates that $\sigma_R$ fixes the elements of $\overline{K_0}$. It contains the element 1 as well as the elements  $a + a^{\sigma_R}$ for all $a \in K_R$. This as $s \in K_0$, and  $b +b^{\sigma s} = s^{-1} (sb + (sb)^\sigma) \in s^{-1} K_0$ for $b \in R$.  

If $a \in s^{-1} K_0 \cap R$ and $b \in R^*$ then $b^{\sigma s} a b \in R$ because $b^{\sigma s} \in R$ by our assumptions. Also $b^{\sigma s} a b  \in s^{-1} K_0$ as $K_0$ is an involutory set and hence contains $b^\sigma sa b$. This implies that if $t\in K_R$, then $t^{\sigma_R} \overline{K_0} t = \overline{K_0}$. Hence $\overline{K_0}$ is an involutory set of $K_R$.
\qed

%When $s^{-1} K_0 \cap R \subset m$ we are in Case II. We have two possibilities for $1+s^{-1}s^\sigma =  s^{-1}(s+s^\sigma) \in s^{-1}K_\sigma \subset s^{-1}K_0$. It can either be an element of $m$, or be an element of $K \setminus R$. Suppose the latter holds. Then $s^{-1}s^\sigma$ also belongs to $k \setminus R$, and the inverse $s^{-\sigma} s$ belongs to $m$. But $(s^{-1}s^\sigma)^{\sigma s} = s^{- \sigma} s$ which contradicts the first compatibility condition. So $1+s^{-1}s^\sigma \in m$ and consequently $s^{-1}s^\sigma \in m$.  Now for a given $a \in R$ we have that $s^{-1}(s^\sigma a + a^\sigma s ) = s^{-1} s^\sigma a + a^{\sigma s} \in s^{-1}K_0 $. This element lies in $m$ as it can not be an element of $k\setminus R$ (by the first compatibility condition and what we already had derived for $s^{-1}s^\sigma$). In particular this implies that $a \equiv a^{\sigma s} \mod m$, or that $\sigma_R$ is the identity automorphism of $K_R$. Because it is also an anti-automorphism this yields that $K_R$ is a field. 

\subsubsection{Structures on $L_0$}

Consider the following two subsets of $L_0$.
\begin{align*}
%L_0' :=& \{v \in L_0 \vert 0< |(q_0(v) -  K_0) \cap sR | \} \\
%L_0'' :=& \{v \in L_0 \vert 0< |(q_0(v) -  K_0) \cap sm | \} 
L_0' :=& \{v \in L_0 \vert  (\exists a \in R) (q_0(v) \equiv sa \mbox{ mod } K_0 )  \} \text{ and} \\
L_0'' :=& \{v \in L_0 \vert (\exists a \in m) (q_0(v) \equiv sa \mbox{ mod } K_0 ) \}. 
\end{align*}

\begin{lemma}
The sets $L_0'$ and $L_0''$ are additive abelian subgroups of $L_0$.
\end{lemma}
\proof
We only proof that $L_0'$ is a subgroup of $L_0$ as the proof for $L_0''$ is completely analogous. Let $v,w \in L_0'$. As $q_0(-v) = q_0(v)$ it is clear that that $L_0'$ is closed under taking inverses. By construction of $L_0'$ we can find $a, b \in sR$ such that $q_0(v) \equiv sa \mbox{ mod } K_0 $ and  $q_0(w) \equiv sb \mbox{ mod } K_0$. By the definition of a skew-hermitian pseudo-quadratic form we have that $$q_0(v+w) \equiv sa +sb +f (v,w) \mbox{ mod } K_0. $$ Condition (C2) asserts that $f(v,w) \in sR$. Hence $sa +sb +f (v,w) \in sR$ and $L_0'$ is indeed a subgroup of $L_0$. \qed

The next lemma investigates how these subgroups behave under scalar products.
\begin{lemma}\label{lemma:mod}
The subgroups $L_0'$ and $L_0''$ are $R$-modules, in particular we have that $L_0' . R = L_0'$ and $L_0'' . R = L_0''$. Moreover  we have that $L_0' . m \subset L_0''$.
\end{lemma}
\proof
Let $v \in L_0'$ and $t \in R$. By construction of $L_0'$ there exists an element $a \in R$ such that $q_0(v)  \equiv sa \mbox{ mod } K_0 $. Then $q_0(v t) \equiv t^\sigma sa t \mbox{ mod } K_0$ as $q_0$ is a skew-hermitian pseudo-quadratic form. From Condition (C1) it follows that $s^{-1} t^\sigma s \in R$, so $ s^{-1}t^\sigma s a t \in R$ or equivalently $t^\sigma sa t \in sR$. This implies that $vt \in L_0'$ and $L_0' . R = L_0'$ in general. The proofs for $L_0'' . R = L_0''$ and $L_0' . m \subset L_0''$ are completely analogous.  \qed

%If one moreover has that $t \in m$, then one easily observes that $t^\sigma sb t \ in sm$ and $vt \in L_0''$. This proves the lemma.

Let $\overline{L_0}$ be the quotient $L_0' / L_0''$. If $v \in L_0' $ and $k \in R$ then $(v + L_0'').(k+m) \subset vk + L_0''$ (by Lemma~\ref{lemma:mod}), so this group can be interpreted as a right vector space over the residue skew field $K_R$.

We construct two functions on $\overline{L_0}$. As first function we define 
$$ \overline{f_0} : \overline{L_0}  \times \overline{L_0} \to K_R:  (v + L_0'',w +L_0'') \mapsto s^{-1} f(v,w) +m. $$
Note that this is indeed a map into $K_R$ by Condition (C2). As second function we define an unary function  $\overline{q_0} : \overline{L_0} \to K_R$, mapping elements $v + L_0''$ of $\overline{L_0}$ to $t + m \in K_R$ (with $t \in R$) such that $q_0(v) \equiv st \mbox{ mod }  K_0$. Note that such an element $t+m$ always exists by definition of $L_0'$, but that $\overline{q_0}$ is not necessarily well-defined as there might be several  $t + m \in K_R$ satisfying this condition and that the possible $t + m \in K_R$ may depend on the choice of representative of $v + L_0''$. However the next three lemmas show that there is indeed a unique choice (albeit modulo $\overline{K_0}$ in the first case), and determine the structure of $\overline{f_0}$ and $\overline{q_0}$.

\begin{lemma}\label{lemma:step1}
Let $v, w\in L_0'$ such that $w+L_0''= v + L_0'' $. If $t \in R$ is such that $q_0(v) \equiv st \mbox{ mod }  K_0$, then there exists a $t' \in t+m$ such that $q_0(w) \equiv st' \mbox{ mod }  K_0$.
\end{lemma}
\proof
Let $v,w, t$ be as in the statement of the lemma, note that $w-v  \in L_0''$. By definition of $L_0''$ there exists an $a \in m$ such that $(w- v, sa) \in T$, or equivalently $q_0(w- v) \equiv sa \mbox{ mod }K_0 $. Note that $f_0(v,w-v) \in sm$ by Condition (C3). As $q_0$ is a pseudo-hermitian form we have
\begin{align*}
q_0(w) &\equiv q_0(v + (w-v)) \mbox{ mod }K_0  \\
&\equiv  q_0(v) + q_0(w-v) + f_0(v,w-v) \mbox{ mod }K_0 \\
&\equiv  st  + sa + f_0(v,w-v) \mbox{ mod }K_0. 
\end{align*}
As $ sa + f_0(v,w-v)$ is in $sm$, the element $t' := t+ a + s^{-1}f_0(v,w-v)$ is in $t+m$ and satisfies $q_0(w) \equiv st' \mbox{ mod }  K_0$. \qed

\begin{lemma} \label{lemma:case1}
Under the assumption of Case I, the map $\overline{f_0}$ is a skew-hermitian sesquilinear function, and $\overline{q_0}$ is well-defined modulo $\overline{K_0}$ and an anisotropic  skew-hermitian pseudo-quadratic form on $\overline{L_0}$ with respect to the involution $\sigma_R$, the involutory set $\overline{K_0}$, and $\overline{f_0}$.
\end{lemma}
\proof
From Lemma~\ref{lemma:step1} and the construction of $\overline{q_0}$ and $\overline{K_0}$, it follows that the function $ \overline{q_0}$ is well-defined modulo $\overline{K_0}$.

The remainder of the statement of the lemma follows from straightforward calculations using the properties of skew-hermitian sesquilinear and pseudo-quadratic forms, Conditions (C1)-(C3) and the fact that $s \in K_0$ (and hence fixed by $\sigma$) by Lemma~\ref{lemma:cases}.
\qed

%  If there are two elements $t+m$ and $t'+m$ in $K_R$ such that for an element $v + L_0''$ of $\overline{L_0}$ such that $st, st' \in q_0(v) -  K_0$, then $t-t' \in s^{-1} K_0$. So $t+m$ and $t'+m$ can only differ by an element of $\overline{K_0}$.

\begin{lemma}
Under the assumption of Case II, we have that $f_0$ is a symmetric bilinear function, and that $\overline{q_0}$ is a well-defined quadratic form on $\overline{L_0}$ with $f_0$ as associated bilinear function.
\end{lemma}
\proof
For an element $v\in L_0'$ let $t,t' \in R$ be two elements such that $q_0(v) \equiv st \equiv st' \mbox{ mod } K_0$. Because $ K_0 \cap sR^*  = \varnothing$ (see Lemma~\ref{lemma:cases}), the difference $st - st' \in K_0$ lies in $sm$, implying that $t +  m= t' +m$ and that, in the light of Lemma~\ref{lemma:step1}, $\overline{q_0}(v + L_0'')$ is well-defined.

Let $v + L_0'',w +L_0''$ be two elements of $\overline{L_0}$. Making use of  Lemma~\ref{lemma:cases} and Condition (C1) we can derive the following.% and that $1+s^{-1}s^\sigma \in m$).
\begin{align*}
\overline{f_0}(v,w) &= \overline{f_0}(v,w)^{\sigma_R} \\
&= (s^{-1} f_0(v,w) +m)^{\sigma s} \\
&= s^{-1} f_0(v,w)^\sigma s^{-\sigma} s + m  \\
&= s^{-1} f_0(v,w)^\sigma + m  \\
&=  s^{-1} f_0(w,v) +m \\
&= \overline{f_0}(w,v)  
\end{align*}
We can conclude that $\overline{f_0}$ is symmetric. The other part of the statement follows easily.
\qed

%Consider the following two subsets of $L_0$.
%\begin{align*}
%L_0' :=& \{v \in L_0 \vert q_0(v) \in sR \}, \\
%L_0'' :=& \{v \in L_0 \vert q_0(v) \in sm \} .
%\end{align*}
%These are both additive abelian groups. Let $\overline{L_0}$ be the quotient $L_0' / L_0''$. If $v \in L_0' $ and $k \in R$ then $(v + L_0'').(k+m) \subset vk + L_0''$, so we this group can be interpreted as a right vector space over $K_R$.
%
%Now define the following two functions:
%\begin{align*}
%\overline{q_0} &: \overline{L_0} \to K_R: v +L_0''\mapsto s^{-1} q(v) + m, \\
%\overline{f_0} &: \overline{L_0}  \times \overline{L_0} : (v + L_0'',w +L_0'') \mapsto s^{-1} f(v,w) +m.
%\end{align*}
%Note that this two functions are well defined by the properties of skew-hermitian sesquilinear and pseudo-quadratic forms. If we are in Case I, then these new maps will be again an anisotropic pseudo-quadratic form and an associated sesquilinear form, and in Case II one obtains an anisotropic quadratic form and an associated symmetric bilinear form.
%

\subsubsection{Constructing the epimorphism }

We will now construct an epimorphism $\rho$ from the polar space $\mathsf{BC}_l(K, K_0, \sigma, L_0, q_0)$ to the polar space $\mathsf{BC}_l(K_R, \overline{K_0}, \sigma_R, \overline{L_0}, \overline{q_0})$ when we are in Case I, and to $\mathsf{B}_l(K_R, \overline{L_0}, \overline{q_0})$ in Case II. We use the notations from Section~\ref{section:psf}.

We call a vector $(v  \vert a_1, \dots, a_{2l}) \in X$ \emph{normed} if all coefficients $a_i$ ($i \in \{1, \dots, 2l\}$) lie in the subring $R$, and at least one is a unit of $R$. The next lemma deals with the existence of  a normed scalar multiple of a given vector. 

\begin{lemma}\label{lemma:norm}
If $w :=(v  \vert a_1, \dots, a_{2l}) \in X$ is a vector such that $( a_1, \dots, a_{2l}) \neq (0, \dots, 0)$, then there exists an element $t \in K$ such that the right scalar product $w t$ is normed. 
\end{lemma}
\proof
Given such a non-zero vector $w := (v  \vert a_1, \dots, a_{2l})$, one can assume without loss of generality (by taking an appropriate right scalar product) that the set $J := \{i \in \{1,\dots, 2l\} \vert a_j \in K \setminus R \}$ is non-empty. Choose a $j \in J$, and let  $w'$ be the vector $w a_{j}^{-1}$. As $a_{j}^{-1} \in m$, this implies that if a coordinate of $w$ was already in $ R$, then this holds for the corresponding coordinate of $w'$ as well. Note that the coordinate of $w'$ corresponding with $j$ is 1, so repeating this algorithm a finite number of times yields the desired scalar multiple.
\qed

The following lemma shows how the different choices of normed scalar multiples are related.

\begin{lemma} \label{lemma:well-defined}
If $w := (v  \vert a_1, \dots, a_{2l}) \in X$ is a vector such that $( a_1, \dots, a_{2l}) \neq (0, \dots, 0)$ and $t, t' \in K$ are elements  such that the right scalar products $w t$ and $w t'$  are normed, then $t^{-1}  t'\in R^*$.  
\end{lemma}
\proof
We will prove this by contradiction. Without loss of generality one may assume that $t^{-1}  t' \in K \setminus R$. By the definition of being normed, there exists a $j \in \{1,\dots , 2l\}$ such that $a_j t \in R^*$. Then $a_jt' = (a_j t)(t^{-1}  t')$ lies in $K \setminus R$, which is impossible for a normed vector.
\qed

Let the vector $(v \vert a_1, a_2, \dots, a_{2l-1},a_{2l})$ represent a point of $\mathsf{BC}_l(K, K_0, \sigma, L_0, q_0)$. Note that as $q_0$ is anisotropic we have that $(a_1, a_2, \dots, a_{2l-1},a_{2l}) \neq (0, \dots, 0)$.  By Lemma~\ref{lemma:norm} we can choose this vector such that  $(v \vert a_1', a_2',  \dots, a_{2l-1}', a_{2l}')$ is normed, with $a_i'$ equal to $a_i$ when $i$ is odd, and equal to $s^{-1}a_i$ when $i$ is even ($i \in \{1,\dots, 2l\}$). We now have that 
\begin{align*}
a_1^\sigma a_2 + \dots + a_{2n-1}^\sigma a_{2n} &= a_1'^\sigma s a_2' + \dots + a_{2n-1}'^\sigma s a_{2n}' \\
&= s  (a_1'^{\sigma s}  a_2' + \dots + a_{2n-1}'^{\sigma s}  a_{2n}' ) \in s R.
\end{align*}

As the 1-dimensional space spanned by the vector is a point of the polar space, we have that $q_0(v) + a_1^\sigma a_2 + \dots + a_{2n-1}^\sigma a_{2l}  \in K_0$. This implies that $v \in L_0'$. In particular we have that $\overline{q_0}(v+ L_0'') = - a_1'^{\sigma s}  a_2' + \dots + a_{2l-1}'^{\sigma s}  a_{2l}'  + m$ (in Case I this is modulo $\overline{K_0}$, see Lemma~\ref{lemma:case1}). 

Hence $\langle (v+ L_0'' \vert a_1'+m, a_2' +m, \dots, a_{2l}' +m) \rangle$ is a point of the polar space $\mathsf{BC}_l(K_R, \overline{K_0}, \sigma_R, \overline{L_0}, \overline{q_0})$ in case I, and a point of $\mathsf{B}_l(K_R, \overline{L_0}, \overline{q_0})$ in Case II. Lemma~\ref{lemma:well-defined} shows that this point does not dependent on the choice of the vector representing the point of $\mathsf{BC}_l(K, K_0, \sigma, L_0, q_0)$. 

We denote the map we defined from the points of the space $\mathsf{BC}_l(K_R, \overline{K_0}, \sigma_R, \overline{L_0}, \overline{q_0})$ to the points of $\mathsf{BC}_l(K_R, \overline{K_0}, \sigma_R, \overline{L_0}, \overline{q_0})$ or $\mathsf{B}_l(K_R, \overline{L_0}, \overline{q_0})$ by $\rho$. We claim that $\rho$ is the desired epimorphism.

\begin{lemma}
The map $\rho$ is surjective.
\end{lemma}
\proof
Let $(v+ L_0'' \vert a_1'+m, a_2' +m, \dots, a_{2l}' +m) \in \overline{L_0} \times K_r^{2l}$ represent a point of the polar space $\mathsf{BC}_l(K_R, \overline{K_0}, \sigma_R, \overline{L_0}, \overline{q_0})$ or $\mathsf{B}_l(K_R, \overline{L_0}, \overline{q_0})$ depending on the case. By definition of these polar spaces and by the definition of $\overline{q_0}$ and $\overline{K_0}$ we have that $s^{-1}q_0(v) + a_1'^{\sigma s} a_2' + \dots + a_{2l-1}'^{\sigma s} a_{2l}' \in s^{-1}K_0 + m$. Without loss of generality we may assume that there is a $j \in \{1, \dots, 2l\}$ such that $a_j' = 1$, this because at least one of the $a_i'$ ($i \in \{1, \dots, 2l\}$) is non-zero modulo $m$ as $\overline{q_0}$ is anisotropic.

Set $a_1 := a_1'$, $a_2 :=  sa_2'$ and so on. Then 
\begin{align*}
q(v \vert a_1, \dots, a_{2l}) &= q_0(v) + a_1^{\sigma} a_2 + \dots + a_{2l-1}^{\sigma } a_{2l} \\
&= s ( s^{-1}q_0(v) + a_1'^{\sigma s} a_2' + \dots + a_{2l-1}'^{\sigma s} a_{2l}') \in K_0 + sm.
\end{align*}  
Choose a $t \in m$ such that $q(v, a_1, \dots, a_{2l})\equiv st \mbox{ mod }K_0$. Let $b_i := a_i$ for all $i \in \{1, \dots, 2l\}$ except for $b_{j - 1} = a_{j-1} -t^{\sigma s^\sigma}$ when $j$ is even and  $b_{j + 1} = a_{j+1} - st$ when $j$ is odd. For each possibility one obtains:
\begin{align*}
q(v, b_1, \dots, b_{2l}) &= q_0(v) + b_1^{\sigma} b_2 + \dots + b_{2l-1}^{\sigma } b_{2l} \\
&= q_0(v) + a_1^{\sigma} a_2 + \dots + a_{2l-1}^{\sigma } a_{2l} - st \\
&= q(v \vert a_1, \dots, a_{2l}) - st \in K_0,
\end{align*}  
by construction of $t$. Hence $\langle (v\vert b_1, \dots, b_{2l}) \rangle$ is a point of $\mathsf{BC}_l(K, K_0, \sigma, L_0, q_0)$, for which one easily verifies that its image under $\rho$ is the point  $\langle(v+ L_0'',a_1'+m, a_2' +m, \dots, a_{2l}' +m) \rangle$. \qed

\begin{lemma}\label{lemma:col}
The map $\rho$ preserves collinearity.
\end{lemma}
\proof
From the construction, the definition of $\overline{f_0}$, and Section~\ref{section:coll}.
\qed

The next series of lemmas proves that a collinearity preserving surjective map induces an epimorphism of the buildings associated to the polar spaces. In order to simplify notations we denote the polar space $\mathsf{BC}_l(K_R, \overline{K_0}, \sigma_R, \overline{L_0}, \overline{q_0})$ by $\Pi$ and the polar space $\mathsf{BC}_l(K_R, \overline{K_0}, \sigma_R, \overline{L_0}, \overline{q_0})$ or $\mathsf{B}_l(K_R, \overline{L_0}, \overline{q_0})$ (depending on the case) by $\Pi'$.

\begin{lemma}\label{lemma:subsp1}
If there exist for a subspace $\pi'$ of $\Pi'$ a subspace $\pi$ of $\Pi$ whose image under $\rho$ is contained in $\pi'$, then for each subspace $\xi'$ of $\pi'$ there exists a subspace $\xi$ of $\pi$ whose image under $\rho$ is contained $\xi'$ and whose co-dimension in $\pi$ is less than or equal to the co-dimension of $\xi'$ in $\pi'$.
\end{lemma}
\proof
We can assume that $\pi'$ is non-empty. We proof this for a subspace $\xi'$ of co-dimension 1 in $\pi'$, the general case then follows by recursion. Embed $\pi'$ in a generator $\chi'$ of $\Pi'$. By Lemma~\ref{lemma:polar2} there exists a point $p'$ of $\Pi'$ such that the intersection of $\pi'$ and the points collinear with $p'$ in $\chi'$ is exactly $\xi'$. Let $p$ be a point in $\Pi$ mapped to $p'$. Let $\xi$ be the subspace of $\pi$ consisting of those points in $\pi$ collinear with $p$, which is of co-dimension 0 or 1 in $\pi$ (see Lemma~\ref{lemma:polar2}). As $\rho$ preserves collinearity it follows that $\xi$ has the desired properties. \qed

\begin{lemma}\label{lemma:subsp2}
If $\pi$ is a subspace of $\Pi$ whose points are mapped by $\rho$ into a same-dimensional subspace $\pi'$ of $\Pi'$, then the points of $\pi$ are mapped surjectively to the points of $\pi'$. In particular the set of points in $\pi$ cannot be mapped into a lower-dimensional subspace of $\Pi'$.
\end{lemma}
\proof 
The first assertion follows from Lemma~\ref{lemma:subsp1} by considering the points of $\pi'$ as 0-dimensional subspaces. The second assertion follows from embedding the lower-dimensional subspace in a subspace of the same dimension as $\pi$ and then applying the first assertion. \qed
 
% \begin{lemma}\label{lemma:line}
%All the points on one line of $\Pi$ cannot be mapped to a single point of $\Pi'$ by $\rho$.
%\end{lemma}
%\proof
%As in a polar space each point is collinear to at least one point on a given line, this would imply by the previous lemmas that every point in the image is collinear with one point, which is impossible. \qed
 
%\begin{lemma}
%If $\pi$ is a subspace of $\Pi$, then its points are mapped surjectively into a same-dimensional subspace of $\Pi'$.
%\end{lemma} 
%\proof
%We only have to prove that $\pi$ is mapped into a subspace $\pi'$ of the same or smaller dimension. Lemma~\ref{lemma:subsp1} implies that $\pi$ and $\pi'$ are of the same dimension, otherwise there is a non-empty subspace of $\pi$ mapped into the empty subspace which is a clear contradiction. Surjectivity then follows from Lemma~\ref{lemma:subsp2}.
%
%As the images of the points of a generator are mutually collinear, they are contained in a subspace. Hence the lemma holds for all generators of $\Pi$. Now assume that $\xi$ is a subspace of $\Pi$ for which lemma holds, and that $\chi$ is a subspace of co-dimension 1 in $\xi$. By embedding $\xi$ in a generator and Lemma ... there exists a point $p$ of $\Pi$ such that the points of $\xi$ collinear with $p$ are exactly the points of $\chi$. The image of $\chi$ is contained in that of $\xi$ and 
%

\begin{lemma}\label{lemma:line2}
The map $\rho$ maps lines of $\Pi$ to subsets of lines of $\Pi'$.
\end{lemma}
\proof
Given a line, let $p_1$ and $p_2$ be two points on this line, not mapped to the same point (this is possible by Lemma~\ref{lemma:subsp2}). Let $p_3$ be a third point on this line, then we have to show that this point is mapped to a point on the line through $p_1^\rho$ and $p_2^\rho$. We may assume that $p_1^\rho \neq p_3^\rho \neq p_2^\rho$. Let $w_1$, $w_2$ and $w_3$ be vectors satisfying the norming condition as in the definition of $\rho$, representing respectively $p_1$, $p_2$ and $p_3$.

Because $p_1$, $p_2$ and $p_3$ lie on a line, there exist non-zero constants $t_1$ and $t_2$ in $K$ such that $w_3 =  w_1  t_1 +  w_2 t_2$. If we can show that $t_1$ and $t_2$ lie in the total subring $R$, we are done (by construction of $\rho$). One can assume without loss of generality that $t_2 t_1^{-1} \in R$. As $ w_3 t_1^{-1} = w_1 + w_2 t_2 t_1^{-1}$, one has that $t_1 \in R$. Otherwise $t_1^{-1}$ would be an element of $m$, implying that $p_1$ and $p_2$ are mapped to the same point by $\rho$. It also follows that $t_2 \in R$, as otherwise it is impossible that $w_3$ satisfies the required norming condition.\qed

\begin{lemma}\label{lemma:subsp3}
The map $\rho$ maps subspaces of $\Pi$ into subsets of same-dimensional subspaces of $\Pi'$.
\end{lemma}
\proof
Let $\pi$ be a subspace of $\Pi$, we prove the lemma by induction on the dimension $t$ of $\pi$. The result is immediate for $t =-1, 0$, and follows from Lemma~\ref{lemma:line2} for $t= 1$. Now suppose the result holds for all subspaces of dimension at most $t-1$. Let $\xi$ be a $(t-1)$-dimensional subspace of $\pi$. The image of $\xi$ is the point set of a $(t-1)$-dimensional subspace $\xi'$ of $\Pi'$ by the induction hypothesis and Lemma~\ref{lemma:subsp2}. Lemma~\ref{lemma:subsp2} moreover implies that $\pi$ is not completely mapped into $\xi'$, hence there exists a point $p$ in $\pi$ mapped outside $\xi'$. As $\rho$ preserves collinearity one has that every point of $\xi'$ is collinear with $p^\rho$. Hence $\xi'$ and $p^\rho$ span a $t$-dimensional subspace $\pi'$ of $\Pi'$. Each point of $\pi$ lies on a line meeting both $p$ and $\xi$, so Lemma~\ref{lemma:line2} yields that this point is mapped to a line meeting $p^\rho$ and $\xi^\rho$, which is hence contained in $\pi'$. This proves the lemma.
\qed

%\proof
%Let $\alpha$ be a subspace of $\mathsf{BC}_l(K, K_0, \sigma, L_0, q_0)$, let $\alpha'$ be a same-dimensional subspace containing the images of the points of $\alpha$. Suppose that $y \in \alpha' $ is not such an image. Let $P'$ be a set of points as described in Lemma~\ref{lemma:polar1} w.r.t. $\alpha'$ and $y$. Choose for each point in $P'$ some point of $\mathsf{BC}_l(K, K_0, \sigma, L_0, q_0)$ mapped to it. The points in the set $P$ one obtains in this way are all collinear with some point $x$ in $\alpha$ by Lemma~\ref{lemma:polar2}. The only possibility for the image of $x$ is $y$ contradicting the way we choose $y$. \qed

\begin{prop}
The map $\rho$ induces an epimorphism between the buildings associated to the polar spaces $\Pi$ and $\Pi'$.
\end{prop}
\proof
The combination of Lemmas~\ref{lemma:subsp2} and~\ref{lemma:subsp3} states that there exists for each subspace $\pi$ of $\Pi$ a unique, same-dimensional subspace $\pi'$ of $\Pi'$ such that the image under $\rho$ of the point set of $\pi$  is exactly the point set of $\pi'$. Hence one can extend $\rho$ to subspaces by setting $\pi^\rho = \pi'$. This extension to subspaces clearly preserves the incidence relation (which is containment) between the subspaces.

In order to have an epimorphism of buildings we only need to prove that this induces a surjective map between the sets of chambers of the buildings associated to $\Pi$ and $\Pi'$, or equivalently that for each maximal flag of the polar space $\Pi'$ there is a maximal flag of $\Pi$ mapped to it by $\rho$. By Lemmas~\ref{lemma:subsp1} and~\ref{lemma:subsp2} it suffices to prove that $\rho$ is surjective on  generators.

Pick a generator $\pi$ in $\Pi$, and consider a generator $\chi$ in $\Pi'$ intersecting $\pi^\rho$ in a subspace of co-dimension one. Let $p'$ be a point of $\chi$ not in the intersection. By surjectivity there exists a point $p$ of $\Pi$ such that $p^\rho = p'$. Note it is impossible that $p$ lies in $\pi$. By Lemma~\ref{lemma:polar} there exists a unique generator $\xi$ containing $p$ and intersecting $\pi$ in a subspace of dimension one less. The only possibillity for the image of $\xi$ is $\chi$. As in a building every two chambers can be connected with a gallery, it follows that repeating this argument yields that we obtain each generator as an image.
 \qed

%\subsubsection{Babycase: Involutory quadrangles}
%We will first consider the case of involutory quadrangles. By~\cite[Section 5.1]{Str:*} one knows that for a given root group $U_i$ there are subgroups $V_i$ (consisting of all the root elations which descend) and $U_i$ (consisting of those who descend to the trivial automorphism). The result in~\cite[Section 6.1]{Str:*}  implies

In order to finish the proof of Theorem~\ref{thm:main2} we only need to check that the total subring and left coset one obtains from applying Theorem~\ref{thm:main1} to the epimorphism $\rho$ are indeed $R$ and $sR^*$. In order to achieve this we need to check which portion of the groups $u(i)$ ($i\in \{1,\dots, l\}$) descends. Recall that these groups were explicitly described in Section~\ref{section:rgl2}.

Let $(v\vert a_1,  \dots, a_{2l})$ be a vector of $X$ representing a point of $\mathsf{BC}_l(K, K_0, \sigma, L_0, q_0)$, satisfying the norming condition appearing in the definition of the epimorphism $\rho$.  So the vector $(v\vert a_1', a_2', \dots)$, with $a_j'$ equal to $a_j$ when $j$ is odd, and equal to $s^{-1}a_j$ when $j$ is even ($j \in \{1,\dots, 2l\}$), is normed. Notice, as before, that this implies that $v \in L_0'$.

We start with the groups $u(i)$ with $i \in \{1,\dots,l-1\}$. We perform the calculations for $i=l-1$, the argument for the other possible values is completely analogous. For $k \in K$ one has
\begin{align*}
(v \vert a_1, \dots, a_{2l})^{y_{l-1}(k)} &= (v \vert a_1', s a_2' , \dots, s a_{2l}')^{y_{l-1}(k)} \\
&= (v \vert a_1' +k a_3',s a_2' , a_3, s a_4 -k^\sigma s a_2 , \dots , s a_{2l}'), \\
&= (v \vert a_1' +k a_3',s a_2' , a_3, s (a_4 -k^{\sigma s} a_2) , \dots , s a_{2l}'). \\
\end{align*}
Note that if $k \in R$ we end up with a vector again satisfying the norming condition appearing in the definition of $\rho$. From this it follows that the automorphism $y_{l-1}(k)$ descends if $k \in R$. If $k \in K \setminus R$ then $y_{l-1}(k)$ maps the vectors $(0\vert 1, 0,0, \dots)$ and $(0 \vert 0, 0, 1, 0, \dots)$ (which are mapped to different points by $\rho$) both to vectors mapped to $(0\vert 1, 0, 0, \dots)$ by $\rho$. Hence  $y_{l-1}(k)$ descends if and only if $k \in R$. In the light of Section~\ref{section:nec} and Lemma~\ref{lemma:n1} this implies that the total subring of $K$ given by Theorem~\ref{thm:main1} is exactly the subring $R$.

For the group $u(l)$ we do a similar calculation. For $(w,t) \in T$ one has 
\begin{align*}
(v \vert a_1, \dots, a_{2l})^{y_l(w,t)} &=  (v \vert a_1', s a_2' ,\dots, s a_{2l}') ^{y_l(w,t)} \\
&= (v+w a_1 '\vert a_1', s a_2' -  t a_1  - f_0(w,v), a_3' , \dots , sa_{2l}'),
\end{align*}
By Condition (C2) we have that if $t \in sR$ (and hence $w \in L_0'$) then this vector again satisfies the norming condition and $y_l(w,t)$ descends. If $t \in K \setminus sR$, then both vectors $(0\vert 1,0,0, \dots)^{y_l(w,t)}$ and $(0 \vert 0,1, 0, \dots)^{y_l(w,t)}$ are mapped to $(0 \vert 0,1, 0, \dots)$ by $\rho$, while $(0\vert 1,0,0, \dots)^\rho$ and $(0\vert 1,0,0, \dots)^\rho$ represent different points. We obtain that $y_l(w,t)$ descends if and only if $t \in sR$. Comparing this with Lemma~\ref{lemma:n2} we see that Theorem~\ref{thm:main1} yields the left coset $sR^*$, as desired.

This concludes the proof of Theorem~\ref{thm:main2}.

\end{document}